\documentclass[12pt]{article}
\usepackage{amsmath, amsthm, amssymb}
\usepackage{amsfonts}
\usepackage{epsfig} 

\usepackage{graphics}
\usepackage{subfigure}
\usepackage{graphicx}
\usepackage{color}
\usepackage{epstopdf} 
\usepackage{float} 
\usepackage{mathrsfs} 
\usepackage{relsize} 

\counterwithin*{equation}{section}
\newcommand{\ds}{\displaystyle}

\newcommand{\bq}{\begin{equation}}
\newcommand{\eq}{\end{equation}}
\newcommand{\bqr}{\begin{eqnarray}}
\newcommand{\eqr}{\end{eqnarray}}
\newcommand{\bqrn}{\begin{eqnarray*}}
\newcommand{\eqrn}{\end{eqnarray*}}

\begin{document} 
\begin{center}	
{\bf\large	 "Phase Transition" in fractional differential equations}\\ 
\text {  }\\
 Pavel B. Dubovski and Jeffrey A. Slepoi\\
 Stevens Institute of Technology, Hoboken, NJ 07030, USA\\
 e-mail: pdubovsk@stevens.edu

\end{center}

\pagenumbering{arabic}
\counterwithout{equation}{section} 
\thispagestyle{empty}
\pagenumbering{roman}
\thispagestyle{plain}
\pagenumbering{arabic}
\addtolength{\parskip}{.1in} 

\noindent{\bf \large Abstract}\\
The basic result in the theory of linear ODEs is that the dimension of the fundamental set of linearly independent solutions is equal to the order of equation. This view is typically extended from integer-order to fractional differential equations. The purpose of this research is to demonstrate that this expectation is incorrect and explain why. We show that linear fractional differential equations may admit additional linearly independent solutions beyond their order. This phenomenon occurs when the coefficients of the equation cross certain threshold values. This “phase transition” produces a hyper-dimensional fundamental set of solutions and offers new insight into fractional differential equations and leads to the revisions of the statements of initial and boundary value problems.\\


\noindent{\bf \large Introduction}
Fractional differential equations are a broader framework for systems whose behavior cannot be properly explained by their instantaneous state alone. For example, fractional models are widely used to describe creep, stress, rheology of materials, (e.g., \cite{2010Meral, 2020Bonfanti}, 
in heterogeneous systems with sub- and superduffusion, \cite{2000Metzler}, 
in control systems with fractional-order controllers, \cite{2013El-Khazali}, 
in electrochemical systems, \cite{2022Lasia}, 
in thermal memory and hereditary stresses in heat-conduction and fluid models, \cite{2015Povstenko}, 
in biology-tissue mechanics, \cite{2006Freed, 2009Coussot}, 
in signal processing and imaging \cite{2016Yang}. 
That's why mathematical investigation of fractional equations is critical for proper analysis and understanding the described complex systems, including chaotic and nonlocal phenomena.      

In this paper we prove that $n^{\rm th}$-order homogeneous Cauchy-Euler fractional equation
\begin{equation}
		\sum_{i=1}^{m}d_i x^{\alpha_i}D_{0\scriptscriptstyle{+}}^{\alpha_i} u(x)=0,\quad x > 0, \label{CEeqn}
\end{equation}
$0\leq n-1<\alpha_i \leq n$, $d_i\in \mathbb{R}$ may have $n+1$ linearly independent solutions. 
 The fractional derivatives are understood in the Riemann-Liouville sense with zero lower terminal value: \cite{2006Kilbas, 2007Boyadjiev}:
 \[
D_{0\scriptscriptstyle{+}}^{\alpha}u(x)=\frac{1}{\Gamma(n-\alpha)}\frac{d^{n}}{dx^n}\int_{0}^{x}\frac{u(t) \, dt}{(x-t)^{\alpha+1-n}}. 
\]

\noindent{\bf \large Characteristic equation}\\ 
We search for the solutions in the form of the fractional power functions
\[
u(x)=x^\gamma, \quad \gamma>-1.
\]
Using formula   
\[
	D_0^\alpha x^\gamma = \frac{\Gamma(1+\gamma)}{\Gamma(1+\gamma-\alpha)}x^{\gamma-\alpha},
\]
we arrive at the following characteristic equation \cite{2021DuSl-3} 
\begin{equation}
	\sum_{i=1}^{m}\frac{d_i\cdot \Gamma(1+\gamma)}{\Gamma(1+\gamma-\alpha_i)}=0. \label{charCE}
\end{equation}

\noindent{\bf \large Roots of the characteristic equation}\\[2mm]
{\bf Lemma 1.} 
Let $n-1\! <  \alpha_i \leq  n$, for all $i=1,...,m$. Then
equation (\ref{charCE}) has at least $n\!-\!1$ roots for $\gamma>-1$.\\[2mm]
{\it Proof.}
\noindent Every addend in (\ref{charCE}) for $\gamma>-1$ has $n$ zeros $\gamma=\{\alpha-1,\ldots,\alpha-n\}$.

\noindent From the Euler reflection formula $
\Gamma(1-z)\Gamma(z)=\ds\frac{\pi}{\sin(\pi z)}, 
$
the characteristic equation (\ref{charCE}) becomes
\begin{equation}
s=\sum_{i=1}^{m} d_i \Gamma(\alpha_i-\gamma) \sin(\pi (\alpha_i-\gamma) )=0. \label{charCE2}
\end{equation}

\begin{figure}[H]
	\begin{center}
		\includegraphics[width=10cm]{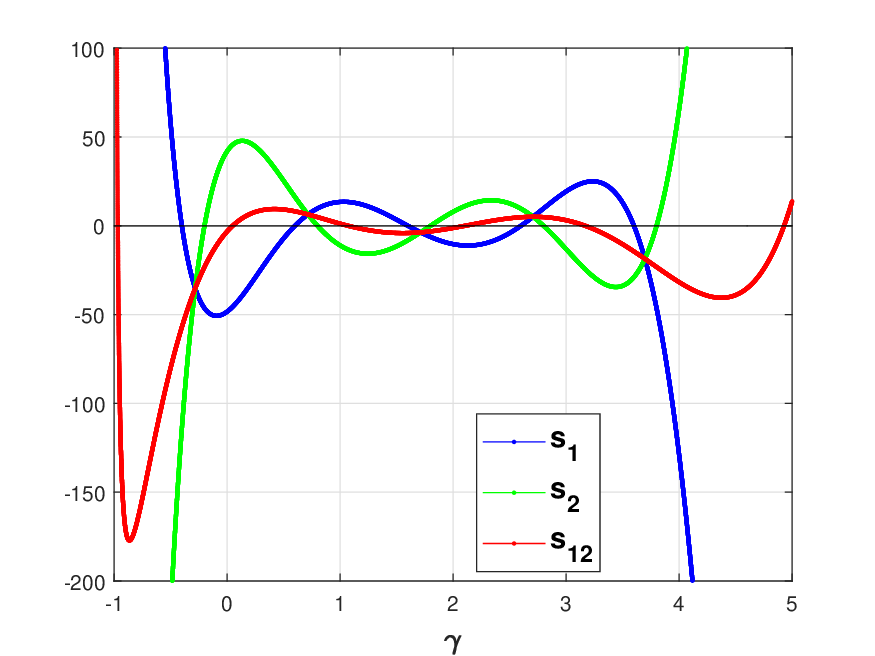}
		\caption{Red curve is the average of graphs $s_1$ and $s_2$ with $d_1=-3.8, d_2=4.0$ and $\alpha_1=4.6, \alpha_2=4.8$.}	\label{FigMultRoots1}
	\end{center}
\end{figure}

%
%
For every addend in (\ref{charCE2}), its $n$ zeros are equidistant along the $\gamma$-axis. Moreover, the zeros of any two addends are shifted by $|\alpha_i-\alpha_j|<1$, keeping the same distance between any pair of neighboring zeros. Thanks to this type of periodicity,  
the intersections of their graphs always change sign. We have one intersection in each of the following intervals $(-1,0],(0,1],\ldots,$ \linebreak $(n-2,n-1]$, and, thus, the graphs of functions
\begin{eqnarray*}
	s_1=d_1\Gamma(\alpha_1-\gamma)\sin(\pi(\alpha_1-\gamma))\ {\rm and}\
	s_2=d_2\Gamma(\alpha_2-\gamma) \sin(\pi(\alpha_2-\gamma))
	\end{eqnarray*}
    have $n$ intersections. 
The graph of their average $\displaystyle s_{12}=\frac{s_1+s_2}{2}$ goes through those $n$ intersection points with alternating sign. Then $s_{12}$ has at least $n-1$ zeros located between $n$ intersections. Figure \ref{FigMultRoots1} clarifies this discussion.  

\vspace{-1cm} \noindent Next, we consider $s_{123}$, the average of $s_1, s_2$, and $s_3$.
Similarly, we observe that $s_3$ and $s_{12}$ have $n$ sign-alternating intersections on $(-1, n-1]$. Hence, $s_{123}=\frac23s_{12}+\frac13 s_3$ also has  at least $n-1$ roots. Repeating this argument for all addends, we see that function $s$ has at least $n-1$ zeros on $(-1,n-1]$. \hfill $\blacksquare$
%

The two unknown intervals remain $(n-1, \infty)$ and the right neighborhood of $\gamma=-1$.  In both cases it is possible that $s$, the left-hand side of equation (\ref{charCE2}), will cross the $ \gamma$-axis (or not cross).  For $\gamma>n-1$, the terms in \eqref{charCE2} with different signs can cause another zero. However, 
if the signs of all coefficients $d_i$ are the same, no additional zeros exist for  $\gamma>n-1$.\\

\noindent{\bf \large Threshold for coefficients}\\[2mm]
{\bf Lemma 2} 
	If the sign of at least one coefficient $d_k$ differs from others, then there exists a threshold coefficient  $D$ such that if $d_k$ crosses $D$, then characteristic equation (\ref{charCE}) has an additional $(n\!+\!1)^{\rm st}$ root close to $\gamma=-1$.      

\begin{figure}[H]
	\begin{center}
		\includegraphics[width=6.5cm]{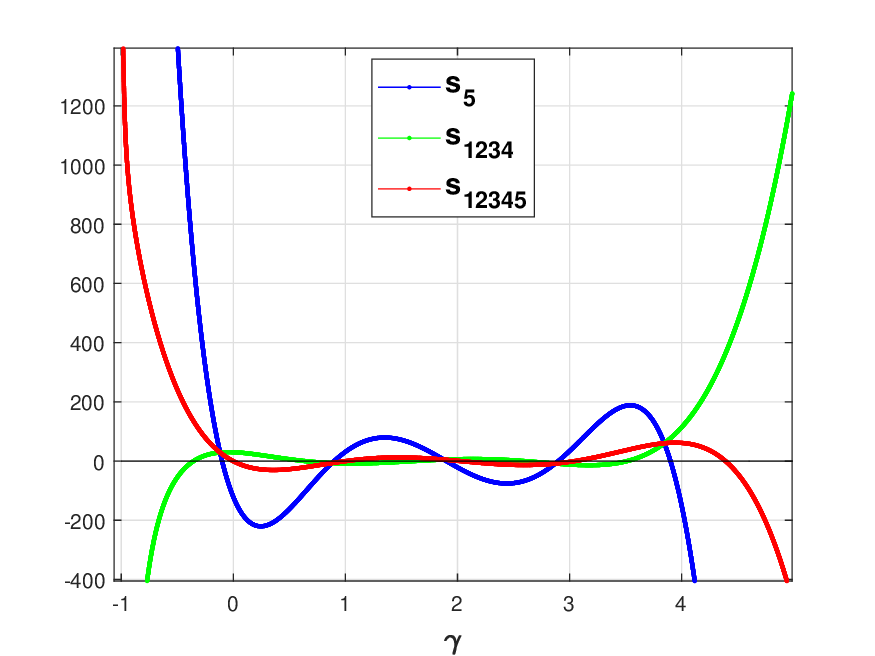}
		\includegraphics[width=6.5cm]{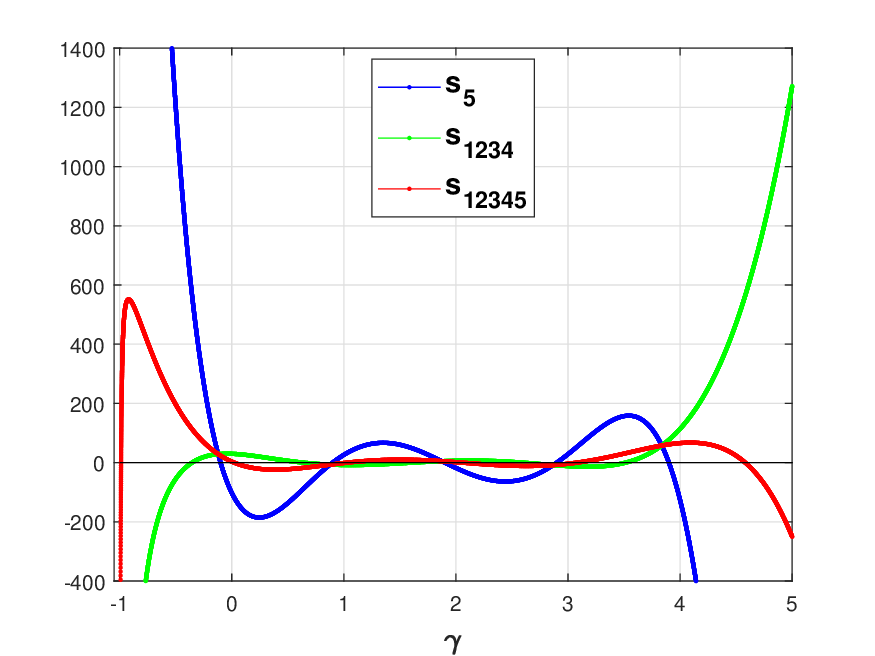}
		\caption{
Left graph: \underline{five} roots (red curve) for $d_5\!=\!-19\!<D\!=\!-17.5761$. \\
Right graph: \underline{six} roots (red curve). Now $d_5\!=\!-16\!>\!D\!=\!-17.5761$. 
		\label{FigMultRoots7}}
	\end{center}
\end{figure}

\noindent{\it Proof.}
We assume that there exists a finite limit $C\neq\infty$ of the left-hand side of  characteristic equation (\ref{charCE}):
\bq
\lim_{\gamma \to -1}\frac{d_k (1+\gamma)}{\Gamma(1+\gamma-\alpha_k)} + \lim_{\gamma \to -1}\sum_{i=1, i\ne k}^{m}\frac{d_i (1+\gamma)}{\Gamma(1+\gamma-\alpha_i)} =C.\nonumber
\eq
Then
\bq
\lim_{\gamma \to -1}\frac{d_k}{\Gamma(1+\gamma-\alpha_k)}=\lim_{\gamma \to -1}\frac{C}{\Gamma(1+\gamma)}-\lim_{\gamma \to -1}\sum_{i=1, i\ne k}^{m}\frac{d_i}{\Gamma(1+\gamma-\alpha_i)}.\nonumber
\eq
Since $\Gamma(0)=\infty$, we obtain $\displaystyle \frac{d_k}{\Gamma(-\alpha_k)}=-\sum_{i=1, i\ne k}^{m}\frac{d_i}{\Gamma(-\alpha_i)}$, and we have arrived at the threshold value for the $k$-th term
\bq
D=-\Gamma(-\alpha_k)\sum_{i=1, i\ne k}^{m}\frac{d_i}{\Gamma(-\alpha_i)}. 
\eq
If we choose $d_k$ higher or lower than the threshold $D$, the limit of the sum in the characteristic equation \eqref{charCE} switches from $-\infty$ to $+\infty$, as can be seen in Figure \ref{FigMultRoots7}. One of those cases generates an extra root near $\gamma=-1$. 

\noindent Below  is the illustrating example for the fifth-order equation
\bq
3.8 D^{4.6}u+4D^{4.8}u + 2D^{4.5}u +2.5 D^{4.1}u +d_5 D^{4.9}u =0. \label{Ex0}
\eq
\vspace{-8mm}\noindent The red curve in the left graph in Figure \ref{FigMultRoots7} shows the case of $n=5$ zeros, the same as the order of equation (\ref{Ex0}). In this case, the value of the fifth coefficient $d_5=-19$ is less than the threshold value $D\!=\!-17.5761$, an extra root is not generated.\\
The red curve in the right graph shows the case of $n=6$ zeros for the same equation (\ref{Ex0}) with a slightly modified coefficient $d_5$. In this case, the value of the fifth coefficient $d_5=-16$ is greater than the threshold value $D\!=\!-17.5761$, which generates the additional root. \hfill $\blacksquare$

\noindent Lemmas 1 and 2 lead us to the following result:\\[2mm]
\noindent{\bf \large Theorem.} 
Let $n-1 < \alpha_i \leq n$, for all $i=1,...,m$. Then, depending on the coefficients $d_i$, equation (\ref{CEeqn}) may have $n+1$ linearly independent solutions.

\noindent {\bf \large Example 1. \underline{Second}-order equation with {\underline{three}} solutions}\\[2mm]
We consider equation
\begin{equation}\label{Example5}
	-8x^{1.6}D_R^{1.6}u(x)+10x^{1.4}D_R^{1.4}u(x)=0.
\end{equation}
Its characteristic equation 
\[
	-\frac{8\Gamma(1+\gamma)}{\Gamma(1+\gamma-1.6)}+\frac{10\Gamma(1+\gamma)}{\Gamma(1+\gamma-1.4)}=0\label{Ex1char}
\]
has three roots: $\gamma_1= -0.9578, \gamma_2=0.1113, \gamma_3=4.0388$. Hence, equation \eqref{Example5} has three independent solutions $u_1=x^{-0.9578}$, $u_2=x^{0.1113}$, $u_3=x^{4.0388}$.\\ 

\noindent{\bf \large Example 2. \underline{First}-order equation with \underline{two}  solutions}\\[2mm]
As another example, we consider first-order equation with integer highest derivative
\begin{equation}\label{Example3}
	xu'(x)-0.5x^{0.5}D^{0.5}u(x)=0.
\end{equation}
Characteristic equation 
\begin{equation*}
	\frac{\Gamma(1+\gamma)}{\Gamma(1+\gamma-1)}-\frac{\Gamma(1+\gamma)}{2\Gamma(1+\gamma-0.5)}=0
\end{equation*}
has two roots $\gamma_1 = -0.86038$ and $\gamma_2=0.42121$ and, thus, this
first-order equation  
    has two linearly independent solutions.\\

\noindent{\bf \large Initial and Boundary value problems}\\[2mm] 
In view of the hyper-dimensionality of the fundamental set, the existence of the extra $(n\!+\!1)^{\rm st}$ solution for the $n^{\rm th}$-order fractional differential equation, we should be very careful in stating the corresponding initial and boundary value problems. $n$ initial or boundary conditions may be insufficient for selecting a unique solution. Particularly, the general solution to (\ref{Example3}) depends on two constants,  $u(x)=c_1x^{-0.86038} +c_2 x^{0.42121}$, and for its well-posedness we should apply an additional condition. Then the boundary value problem to the first-order equation 
\[
xu' - \frac1{2}\sqrt{x} D^{1/2} u=0,\quad u(x_1)=a,\ u(x_2)=b,
\]
becomes well posed for any boundary conditions $a,\ b$ and positive endpoints $x_1,\ x_2$.\\

\noindent {\bf \large Conclusions}\\
We prove that the dimension of the fundamental set of linearly independent solutions for linear fractional differential equations may be higher than the order of the equation if its coefficients pass certain numerical thresholds. This hyper-dimensionality phenomenon drastically influences how to define the corresponding initial/boundary value problems and Green functions. To provide well-posed statements of the problems, the number of initial or boundary conditions must depend not only on the order of equations but also on the numerical values of the coefficients.


\nocite{*}

\end{document}